\documentstyle[12pt]{article}
\begin{document}

\newcommand{\mysection }[1]{\section{#1}\setcounter{equation}{0}}

\title{\bf A Simple Proof for the Generalized Frankel Conjecture} \vspace{15mm}
\author{\bf}
\date{}
\maketitle \centerline{\large \bf Hui-Ling Gu } \vspace{8mm}
\centerline{{\large Department of Mathematics}} \vspace{3mm}
\centerline{{\large Sun Yat-Sen University }} \vspace{3mm}
\centerline{{\large Guangzhou, P.R.China}} \vspace{3mm}
\vspace{15mm}

\noindent {\bf Abstract } In this short paper, we will give a
simple and transcendental  proof for Mok's theorem of the
generalized Frankel conjecture. This work is based on the maximum
principle in \cite{BS2} proposed by Brendle and Schoen.

\begin{center}

{\bf \large \bf{1. Introduction}}
\end{center}

Let $M^n$ be an $n$-dimensional compact K$\ddot{a}$hler manifold.
The famous Frankel conjecture states that: if $M$ has positive
holomorphic bisectional curvature, then it is biholomorphic to the
complex projective space $CP^n$. This was independently proved by
Mori \cite{Mori} in 1979 and Siu-Yau \cite{Siu-Yau} in 1980 by using
different methods. Mori had got a more general result. His method is
to study the deformation of a morphism from $CP^1$ into the
projective manifold $M^n$, while Siu-Yau used the existence result
of minimal energy 2-spheres to prove the Frankel conjecture. After
the work of Mori and Siu-Yau, it is natural to ask the question for
the semi-positive case: what the manifold is if the holomorphic
bisectional curvature is nonnegative. This is often called the
generalized Frankel conjecture and was proved by Mok \cite{Mok}. The
exact statement is as follows:

\vskip 0.3cm \noindent {\bf Theorem 1.1} \emph{ Let $(M,h)$ be an
$n$-dimensional compact K$\ddot{a}$hler manifold of nonnegative
holomorphic bisectional curvature and let $(\tilde{M},\tilde{h})$
be its universal covering space. Then there exist nonnegative
integers $k,N_1,\cdots,N_l$ and irreducible compact Hermitian
symmetric spaces $M_1,\cdots,M_p$ of rank $\geq 2$ such that
$(\tilde{M},\tilde{h})$ is isometrically biholomorphic to
$$(C^k,g_0)\times(CP^{N_1},\theta_1)\times\cdots\times(CP^{N_l},\theta_l)
\times(M_1,g_1)\times\cdots\times(M_p,g_p)$$ where $g_0$ denotes the
Euclidean metric on $C^k$, $g_1,\cdots, g_p$ are canonical metrics
on $M_1,\cdots, M_p$, and $\theta_i,1 \leq i \leq l$, is a
K$\ddot{a}$hler metric on $CP^{N_i}$ carrying nonnegative
holomorphic bisectional curvature.}

We point out that the three dimensional case of this result was
obtained by Bando \cite{Ban}. In the special case, for all
dimensions, when the curvature operator of $M$ is assumed to be
nonnegative, the above result was proved by Cao and Chow
\cite{CC}.

By using the splitting theorem of Howard-Smyth-Wu \cite{HSW}, one
can reduce Theorem 1.1 to the proof of the following theorem:

\vskip 0.3cm \noindent {\bf Theorem 1.2} \emph{ Let $(M,h)$ be an
$n$-dimensional compact simply connected K$\ddot{a}$hler manifold
of nonnegative holomorphic bisectional curvature such that the
Ricci curvature is positive at one point. Suppose the second Betti
number $b_2(M)=1$. Then either $M$ is biholomorphic to the complex
projective space or $(M,h)$ is isometrically biholomorphic to an
irreducible compact Hermitian symmetric manifold of rank $\geq
2$.}

In \cite{Mok}, Mok proved Theorem 1.2 and hence the generalized
Frankel conjecture. His method depended on Mori's theory of
rational curves on Fano manifolds, so it was not completely
transcendental in nature. The purpose of this paper is to give a
completely transcendental proof of Theorem 1.2.

Our method is inspired by the recent breakthroughs in Ricci flow
due to \cite{BW,BS1,BS2}. In \cite{BW}, by developing a new method
constructing the invariant cones to Ricci flow,  B$\ddot{o}$hm and
Wilking proved the differentiable sphere theorem for manifolds
with positive curvature operator.  Recently, Brendle and Schoen
\cite{BS1} proved the $\frac{1}{4}$-differentiable sphere theorem
by using method of \cite{BW}.  Moreover in \cite{BS2}, the authors
gave a complete classification of weakly $\frac{1}{4}$-pinched
manifolds. In this paper, we will use  the powerful strong maximum
principle proposed in \cite{BS2} to give Theorem 1.2 a simple
proof.

\vskip 0.3cm \noindent {\bf Acknowledgement}  I would like to
thank my advisor Professor X.P.Zhu and Professor B.L.Chen for
their encouragement, suggestions and discussions. This paper was
done under their advice.

\begin{center}{ \bf \large \bf {2. The Proof of the Main Theorem}}\end{center}
 \vskip 0.1cm \noindent{\bf Proof of the Main Theorem 1.2.} Suppose
 $(M,h)$ is a compact simply connected K$\ddot{a}$hler manifold of
 nonnegative holomorphic bisectional curvature such that the Ricci
 curvature is positive at one point. We evolve the metric by the K$\ddot{a}$hler Ricci flow:
$$
      \left\{
       \begin{array}{lll}
\frac{\partial}{\partial t}g_{i\bar{j}}(x,t)=-R_{i\bar{j}}(x,t),
          \\[4mm]
  g_{i\bar{j}}(x,0)=h_{i\bar{j}}(x).
       \end{array}
    \right.
$$
According to Bando \cite{Ban}, we know that the evolved metric
$g_{i\bar{j}}(t), t\in (0,T)$, remains K$\ddot{a}$hler. Then by
Proposition 1.1 in \cite{Mok}, we know that for $t\in (0,T)$,
$g_{i\bar{j}}(t)$ has nonnegative holomorphic bisectional
curvature and positive holomorphic sectional curvature and
positive Ricci curvature everywhere. Moreover, according to
Hamilton \cite{Ha86}, under the evolving orthonormal frame
$\{e_\alpha\}$, we have
$$\frac{\partial}{\partial t}R_{\alpha\bar{\alpha}\beta\bar{\beta}}=
\triangle R_{\alpha\bar{\alpha}\beta\bar{\beta}}+\Sigma
_{\mu,\nu}(R_{\alpha\bar{\alpha}\mu\bar{\nu}}R_{\nu\bar{\mu}\beta\bar{\beta}}
-|R_{\alpha\bar{\mu}\beta\bar{\nu}}|^2+|R_{\alpha\bar{\beta}\mu\bar{\nu}}|^2).$$

Suppose $(M,h)$ is not locally symmetric. In the following, we
want to show that $M$ is biholomorphic to the complex projective
space $CP^n$.

Since the smooth limit of locally symmetric space is also locally
symmetric, we can obtain that there exists $\delta \in (0,T)$ such
that $(M,g_{i\bar{j}}(t))$ is not locally symmetric for $t \in
(0,\delta)$. Combining the K$\ddot{a}$hlerity of $g_{i\bar{j}}
(t)$ and Berger's holonomy theorem, we know that the holonomy
group Hol$(g(t))=U(n)$.

Let $P=\bigcup_{p\in M}(T_p^{1,0}(M)\times T_p^{1,0}(M))$ be the
fiber bundle with the fixed metric $h$ and the fiber over $p \in
M$ consists of all 2-vectors $\{X,Y\}\subset T_p^{1,0}(M)$. We
define a function $u$ on $P\times(0,\delta)$ by
$$u(\{X,Y\},t)=R(X,\overline{X},Y,\overline{Y}),$$ where $R$
denotes the pull-back of the curvature tensor of $g_{i\bar{j}}(t)$.
Clearly we have $u\geq 0$, since $(M,g_{i\bar{j}}(t))$ has
nonnegative holomorphic bisectional curvature. Denote
$F=\{(\{X,Y\},t)|u(\{X,Y\},t)=0,X\neq 0, Y\neq 0\}\subset P\times
(0,\delta)$ of all pairs $(\{X,Y\},t)$ such that $\{X,Y\}$ has zero
holomorphic bisectional curvature with respect to $g_{i\bar{j}}(t)$.
Following Mok \cite{Mok}, we consider the Hermitian form
$H_\alpha(X,Y)=R(e_\alpha,\overline{e_\alpha},X,\overline{Y})$, for
all $X,Y\in T_p^{1,0}(M)$ and all $p\in M$, attached to $e_\alpha$.
Let $\{E_\mu\}$ be an orthonormal basis associated to eigenvectors
of $H_\alpha$. In the basis we have
$$\sum_{\mu,\nu}R_{\alpha\bar{\alpha}\mu\bar{\nu}}R_{\nu\bar{\mu}\beta\bar{\beta}}=\sum_\mu
R(e_\alpha,\overline{e_\alpha},E_\mu,\overline{E_\mu})R(E_\mu,\overline{E_\mu},e_\beta,\overline{e_\beta}),$$
and
$$\sum_{\mu,\nu}|R_{\alpha\bar{\mu}\beta\bar{\nu}}|^2
=\sum_{\mu,\nu}|R(e_\alpha,\overline{E_\mu},e_\beta,\overline{E_\nu})|^2.$$

First, we claim that:
$$\sum_{\mu,\nu}
R_{\alpha\bar{\alpha}\mu\bar{\nu}}R_{\nu\bar{\mu}\beta\bar{\beta}}-
\sum_{\mu,\nu}|R_{\alpha\bar{\mu}\beta\bar{\nu}}|^2\geq c_1\cdot
\min\{0,\inf_{|\xi|=1,\xi\in
V}D^2u(\{e_{\alpha},e_{\beta}\},t)(\xi,\xi)\},$$ for some constant
$c_1>0$, where $V$ denotes the vertical subspaces.

Indeed, inspired by Mok \cite{Mok}, for any given $\varepsilon_0>0$
and each fixed $\chi\in \{1,2,\cdot\cdot\cdot,n\}$, we consider the
function
$$\widetilde{G}_\chi(\varepsilon)=(R+\varepsilon_0R_0)(e_\alpha+\varepsilon E_\chi,
\overline{e_\alpha+\varepsilon E_\chi},e_\beta+ \varepsilon\sum_\mu
C_\mu E_\mu,\overline{e_\beta+ \varepsilon\sum_\mu C_\mu E_\mu}),$$
where $R_0$ is a curvature operator defined by
$(R_0)_{i\bar{j}k\bar{l}}=g_{i\bar{j}}g_{k\bar{l}}+g_{i\bar{l}}g_{k\bar{j}}$
and $C_\mu$ are complex constants to be determined later. For the
simplicity, we denote $\widetilde{R}=R+\varepsilon_0R_0$, then
$$\widetilde{G}_\chi(\varepsilon)=\widetilde{R}(e_\alpha+\varepsilon E_\chi,
\overline{e_\alpha+\varepsilon E_\chi},e_\beta+ \varepsilon\sum_\mu
C_\mu E_\mu,\overline{e_\beta+ \varepsilon\sum_\mu C_\mu E_\mu}).$$
Then a direct computation gives
$$\arraycolsep=1.5pt\begin{array}{rcl}
&&\frac{1}{2}\cdot\frac{d^2\widetilde{G}_\chi(\varepsilon)}{d\varepsilon^2}|_{\varepsilon=0}
=\widetilde{R}(E_\chi,\overline{E_\chi},e_\beta,\overline{e_\beta})
+\sum_\mu
|C_\mu|^2\widetilde{R}(e_\alpha,\overline{e_\alpha},E_\mu,\overline{E_\mu})\\[4mm]
&&\hskip 2.5cm+2Re\sum_\mu \overline{C_\mu}
\widetilde{R}(e_\alpha,\overline{E_\chi},e_\beta,\overline{E_\mu})+2Re\sum_\mu
C_\mu
\widetilde{R}(e_\alpha,\overline{e_\beta},E_\mu,\overline{E_\chi}).
\end{array}$$
Writing $C_\mu
=x_\mu e^{i\theta_\mu}$, $(\mu \geq 1)$, for $x_\mu, \theta_\mu$ are
constants to be determined later, the above identity is:
$$\arraycolsep=1.5pt\begin{array}{rcl}
&&\frac{1}{2}\cdot\frac{d^2\widetilde{G}_\chi(\varepsilon)}{d\varepsilon^2}|_{\varepsilon=0}
=\widetilde{R}(E_\chi,\overline{E_\chi},e_\beta,\overline{e_\beta})
+\sum_\mu
|x_\mu|^2\widetilde{R}(e_\alpha,\overline{e_\alpha},E_\mu,\overline{E_\mu})\\[4mm]
&&\hskip 2.5cm+2\sum_\mu x_\mu\cdot Re(e^{-i\theta_\mu}
\widetilde{R}(e_\alpha,\overline{E_\chi},e_\beta,\overline{E_\mu})
+e^{i\theta_\mu}\widetilde{R}(e_\alpha,\overline{e_\beta},E_\mu,\overline{E_\chi})).
\end{array}$$
Following Mok \cite{Mok}, by setting $A_\mu
=\widetilde{R}(e_\alpha,\overline{e_\beta},E_\mu,\overline{E_\chi}),
B_\mu
=\widetilde{R}(e_\alpha,\overline{E_\chi},e_\beta,\overline{E_\mu})$,
we have:
$$\arraycolsep=1.5pt\begin{array}{rcl}
&&\frac{1}{2}\cdot\frac{d^2\widetilde{G}_\chi(\varepsilon)}{d\varepsilon^2}|_{\varepsilon=0}
=\widetilde{R}(E_\chi,\overline{E_\chi},e_\beta,\overline{e_\beta})
+\sum_\mu
|x_\mu|^2\widetilde{R}(e_\alpha,\overline{e_\alpha},E_\mu,\overline{E_\mu})\\[4mm]
&&\hskip 2.5cm+\sum_\mu x_\mu(e^{-i\theta_\mu}
B_\mu+e^{i\theta_\mu}\overline{B_\mu}+e^{i\theta_\mu}A_\mu+e^{-i\theta_\mu}\overline{A_\mu})\\[4mm]
&&\hskip
2.5cm=\widetilde{R}(E_\chi,\overline{E_\chi},e_\beta,\overline{e_\beta})
+\sum_\mu
|x_\mu|^2\widetilde{R}(e_\alpha,\overline{e_\alpha},E_\mu,\overline{E_\mu})
\\[4mm]
&&\hskip 2.5cm+\sum_\mu
x_\mu\cdot(\overline{e^{i\theta_\mu}(A_\mu+\overline{B_\mu})}+e^{i\theta_\mu}(A_\mu+\overline{B_\mu}))
\end{array}$$
By choosing $\theta_\mu$ such that
$e^{i\theta_\mu}(A_\mu+\overline{B_\mu})$ is real and positive,
the identity becomes:
$$\arraycolsep=1.5pt\begin{array}{rcl}
&&\frac{1}{2}\cdot\frac{d^2\widetilde{G}_\chi(\varepsilon)}{d\varepsilon^2}|_{\varepsilon=0}
=\widetilde{R}(E_\chi,\overline{E_\chi},e_\beta,\overline{e_\beta})
+\sum_\mu
|x_\mu|^2\widetilde{R}(e_\alpha,\overline{e_\alpha},E_\mu,\overline{E_\mu})\\[4mm]
&&\hskip 2.5cm+2\sum_\mu x_\mu\cdot |A_\mu+\overline{B_\mu}|.
\end{array}$$
If we change $e_\alpha$ with $e^{i\varphi}e_\alpha$, then $A_\mu
=\widetilde{R}(e_\alpha,\overline{e_\beta},E_\mu,\overline{E_\chi})$
is replaced by $e^{i\varphi}A_\mu$, and
$\overline{B_\mu}=\overline{\widetilde{R}(e_\alpha,\overline{E_\chi},e_\beta,\overline{E_\mu})}$
is replaced by $e^{-i\varphi}\overline{B_\mu}$, we have:
$$\arraycolsep=1.5pt\begin{array}{rcl}
&&\frac{1}{2}\cdot\frac{d^2\widetilde{F}_\chi(\varepsilon)}{d\varepsilon^2}|_{\varepsilon=0}
=\widetilde{R}(E_\chi,\overline{E_\chi},e_\beta,\overline{e_\beta})
+\sum_\mu
|x_\mu|^2\widetilde{R}(e_\alpha,\overline{e_\alpha},E_\mu,\overline{E_\mu})\\[4mm]
&&\hskip 2.5cm+2\sum_\mu x_\mu\cdot
|e^{i\varphi}A_\mu+e^{-i\varphi}\overline{B_\mu}|,
\end{array}$$
where
$$\widetilde{F}_\chi(\varepsilon)=\widetilde{R}(e^{i\varphi}e_\alpha+\varepsilon E_\chi,
\overline{e^{i\varphi}e_\alpha+\varepsilon E_\chi},e_\beta+
\varepsilon\sum_\mu C_\mu E_\mu,\overline{e_\beta+
\varepsilon\sum_\mu C_\mu E_\mu}).$$ Since the curvature operators
$R$ and $R_0$ have nonnegative and positive holomorphic bisectional
curvature respectively, we know that the operator
$\widetilde{R}=R_0+\varepsilon_0R_0$ has positive holomorphic
bisectional curvature. Now by choosing
$x_\mu=-\frac{|e^{i\varphi}A_\mu+e^{-i\varphi}\overline{B_\mu}|}{\widetilde{R}(e_\alpha,\overline{e_\alpha},E_\mu,\overline{E_\mu})}$,
for $ \mu\geq 1 $, it follows that
$$\frac{1}{2\pi}\int_0^{2\pi}(\frac{1}{2}\cdot\frac{d^2\widetilde{F}_\chi(\varepsilon)}{d\varepsilon^2}|_{\varepsilon=0})d\varphi
=\widetilde{R}(E_\chi,\overline{E_\chi},e_\beta,\overline{e_\beta})
-\sum_\mu
\frac{|A_\mu|^2+|B_\mu|^2}{\widetilde{R}(e_\alpha,\overline{e_\alpha},E_\mu,\overline{E_\mu})}$$
and then
$$\arraycolsep=1.5pt\begin{array}{rcl}
&&\hskip
0.5cm\widetilde{R}(e_\alpha,\overline{e_\alpha},E_\chi,\overline{E_\chi})\cdot
\frac{1}{2\pi}\int_0^{2\pi}
(\frac{1}{2}\cdot\frac{d^2\widetilde{F}_\chi(\varepsilon)}{d\varepsilon^2}|_{\varepsilon=0})d\varphi\\[4mm]
&&\hskip
0.1cm=\widetilde{R}(e_\alpha,\overline{e_\alpha},E_\chi,\overline{E_\chi})\widetilde{R}(E_\chi,\overline{E_\chi},e_\beta,\overline{e_\beta})
-\sum_{\mu}\frac{|A_\mu|^2+|B_\mu|^2}{\widetilde{R}(e_\alpha,\overline{e_\alpha},E_\mu,\overline{E_\mu})}
\widetilde{R}(e_\alpha,\overline{e_\alpha},E_\chi,\overline{E_\chi}).
\end{array}$$
Note that
$$\arraycolsep=1.5pt\begin{array}{rcl}
&&\widetilde{F}_\chi(\varepsilon)=\widetilde{R}(e^{i\varphi}e_\alpha+\varepsilon
E_\chi, \overline{e^{i\varphi}e_\alpha+\varepsilon E_\chi},e_\beta+
\varepsilon\sum_\mu C_\mu E_\mu,\overline{e_\beta+
\varepsilon\sum_\mu C_\mu E_\mu})\\[4mm]
&&\hskip 1.1cm=\widetilde{R}(e_\alpha+\varepsilon e^{-i\varphi}
E_\chi, \overline{e_\alpha+\varepsilon e^{-i\varphi}E_\chi},e_\beta+
\varepsilon\sum_\mu C_\mu E_\mu,\overline{e_\beta+
\varepsilon\sum_\mu C_\mu E_\mu}).
\end{array}$$
Interchanging the roles of $E_\chi$ and $E_\mu$, and then taking
summation, we have
$$\arraycolsep=1.5pt\begin{array}{rcl}
&&\hskip
0.5cm\sum_{\chi}2\widetilde{R}(e_\alpha,\overline{e_\alpha},E_\chi,\overline{E_\chi})\widetilde{R}(E_\chi,\overline{E_\chi},e_\beta,\overline{e_\beta})
\\[4mm]
&&\hskip 0.1cm\geq c_1\cdot \min\{0,\inf_{|\xi|=1,\xi\in V}D^2\widetilde{u}(\{e_{\alpha},e_{\beta}\},t)(\xi,\xi)\}\\[4mm]
&&\hskip 0.5cm+\sum_{\mu
,\chi}(|A_\mu|^2+|B_\mu|^2)(\frac{\widetilde{R}(e_\alpha,\overline{e_\alpha},E_\chi,\overline{E_\chi})}{R(e_\alpha,\overline{e_\alpha},E_\mu,\overline{E_\mu})}+
\frac{\widetilde{R}(e_\alpha,\overline{e_\alpha},E_\mu,\overline{E_\mu})}{\widetilde{R}(e_\alpha,\overline{e_\alpha},E_\chi,\overline{E_\chi})})\\[4mm]
&&\hskip 0.1cm\geq c_1\cdot \min\{0,\inf_{|\xi|=1,\xi\in
V}D^2\widetilde{u}(\{e_{\alpha},e_{\beta}\},t)(\xi,\xi)\}+2\sum_{\mu
,\chi}|\widetilde{R}(e_\alpha,\overline{E_\chi},e_\beta,\overline{E_\mu})|^2,
\end{array}$$
where
$\widetilde{u}(\{X,Y\},t)=\widetilde{R}(X,\overline{X},Y,\overline{Y})=
R(X,\overline{X},Y,\overline{Y})+\varepsilon_0R_0(X,\overline{X},Y,\overline{Y})$
and $c_1$ is a positive constant that does not depend on
$\varepsilon_0$.

Hence
$$\arraycolsep=1.5pt\begin{array}{rcl}
&&\hskip 0.5cm\sum_{\mu}
\widetilde{R}(e_\alpha,\overline{e_\alpha},E_\mu,\overline{E_\mu})\widetilde{R}(E_\mu,\overline{E_\mu},e_\beta,\overline{e_\beta})
-\sum_{\mu , \nu }
|\widetilde{R}(e_\alpha,\overline{E_\mu},e_\beta,\overline{E_\nu})|^2\\[4mm]
&&\hskip 0.1cm\geq c_1\cdot \min\{0,\inf_{|\xi|=1,\xi\in
V}D^2\widetilde{u}(\{e_{\alpha},e_{\beta}\},t)(\xi,\xi)\}.
\end{array}$$
Since $\varepsilon_0>0$ is arbitrary, we can let $\varepsilon_0
\rightarrow 0$, then we obtain that:
$$\sum_{\mu,\nu}
R_{\alpha\bar{\alpha}\mu\bar{\nu}}R_{\nu\bar{\mu}\beta\bar{\beta}}-
\sum_{\mu,\nu}|R_{\alpha\bar{\mu}\beta\bar{\nu}}|^2\geq c_1\cdot
\min\{0,\inf_{|\xi|=1,\xi\in
V}D^2u(\{e_{\alpha},e_{\beta}\},t)(\xi,\xi)\},$$ for some constant
$c_1>0$. Therefore we proved our first claim.

By the definition of $u$ and the evolution equation of the
holomorphic bisectional curvature, we know that
$$\arraycolsep=1.5pt\begin{array}{rcl}
&&\frac{\partial }{\partial t}u(\{X,Y\},t)=\triangle u(\{X,Y\},t)+
\sum_{\mu,\nu}R(X,\overline{X},e_\mu,\overline{e_\nu})R(e_\nu,\overline{e_\mu},Y,\overline{Y})\\[4mm]
&&\hskip 3cm-\sum_{\mu
,\nu}|R(X,\overline{e_\mu},Y,\overline{e_\nu})|^2+\sum_{\mu
,\nu}|R(X,\overline{Y},e_\mu,\overline{e_\nu})|^2.
\end{array}$$
Combining the above inequality, we obtain that:
$$\frac{\partial u}{\partial t}\geq L u+ c_1\cdot \min\{0,\inf_{|\xi|=1,\xi\in V}D^2u(\xi,\xi)\},$$
where $L$ is the horizontal Laplacian on $P$, $V$ denotes the
vertical subspaces. By Proposition 2 in \cite{BS2}, (Actually, the
same argument still holds for the bundle $P$ in \cite{BS2} changed
by the bundle $P$ defined in our paper.), we know that the set
$$F=\{(\{X,Y\},t)| u(\{X,Y\},t)=0,X\neq 0, Y\neq 0\}\subset P\times
(0,\delta)$$ is invariant under parallel transport.

Next, we claim that $R_{\alpha\bar{\alpha}\beta\bar{\beta}}>0$ for
all $t\in (0,\delta).$

Indeed, suppose not. Then $R_{\alpha\bar{\alpha}\beta\bar{\beta}}=0$
for some $t\in (0,\delta)$. Therefore
$$(\{e_\alpha,e_\beta\},t)\in F.$$
Combining $R_{\alpha\bar{\alpha}\beta\bar{\beta}}=0$ and the
evolution equation of the curvature operator and the first
variation, we can obtain that:
$$
      \left\{
       \begin{array}{lll}
\sum_{\mu,\nu}(R_{\alpha\bar{\alpha}\mu\bar{\nu}}R_{\nu\bar{\mu}\beta\bar{\beta}}
-|R_{\alpha\bar{\mu}\beta\bar{\nu}}|^2)=0,
          \\[4mm]
R_{\alpha\bar{\beta}\mu\bar{\nu}}=0, \quad \forall \mu, \nu,
          \\[4mm]
R_{\alpha\bar{\alpha}\mu\bar{\beta}}=
R_{\beta\bar{\beta}\mu\bar{\alpha}}=0,  \quad \forall \mu.
         \end{array}
    \right.
$$
We define an orthonormal 2-frames
$\{\widetilde{e_\alpha},\widetilde{e_\beta}\}\subset T_p^{1,0}(M)$
by
$$\widetilde{e_\alpha}=\sin \theta\cdot e_\alpha-\cos \theta\cdot e_\beta,$$
$$\widetilde{e_\beta}=\cos \theta\cdot  e_\alpha+\sin \theta\cdot e_\beta.$$
Then
$$\overline{\widetilde{e_\alpha}}=
\sin \theta\cdot \overline{e_\alpha}-\cos \theta\cdot
\overline{e_\beta},$$
$$\overline{\widetilde{e_\beta}}=
\cos \theta\cdot \overline{e_\alpha}+\sin \theta\cdot
\overline{e_\beta}.$$ Since $F$ is invariant under parallel
transport and $(M,g_{i\bar{j}}(t))$ has holonomy group $U(n)$, we
obtain that
$$(\{\widetilde{e_\alpha},\widetilde{e_\beta}\},t)\in F,$$
that is,
$$R(\widetilde{e_\alpha},\overline{\widetilde{e_\alpha}},
\widetilde{e_\beta},\overline{\widetilde{e_\beta}})=0.$$ On the
other hand,
$$\arraycolsep=1.5pt\begin{array}{rcl}
&&R(\widetilde{e_\alpha},\overline{\widetilde{e_\alpha}},
\widetilde{e_\beta},\overline{\widetilde{e_\beta}})=\sin^2\theta\cos^2\theta
R_{\alpha\bar{\alpha}\alpha\bar{\alpha}}+\sin^3\theta\cos\theta
R_{\alpha\bar{\alpha}\alpha\bar{\beta}}+\sin^3\theta\cos\theta
R_{\alpha\bar{\alpha}\beta\bar{\alpha}}\\[4mm]
&&\hskip 3.3cm+\sin^4\theta
R_{\alpha\bar{\alpha}\beta\bar{\beta}}-\sin\theta\cos^3\theta
R_{\alpha\bar{\beta}\alpha\bar{\alpha}}-\sin^2\theta\cos^2\theta
R_{\alpha\bar{\beta}\alpha\bar{\beta}}\\[4mm]
&&\hskip 3.3cm-\sin^2\theta\cos^2\theta
R_{\alpha\bar{\beta}\beta\bar{\alpha}}-\sin^3\theta\cos\theta
R_{\alpha\bar{\beta}\beta\bar{\beta}}-\cos^3\theta\sin\theta
R_{\beta\bar{\alpha}\alpha\bar{\alpha}}\\[4mm]
&&\hskip 3.3cm-\sin^2\theta\cos^2\theta
R_{\beta\bar{\alpha}\alpha\bar{\beta}}-\sin^2\theta\cos^2\theta
R_{\beta\bar{\alpha}\beta\bar{\alpha}}-\cos\theta\sin^3\theta
R_{\beta\bar{\alpha}\beta\bar{\beta}}\\[4mm]
&&\hskip 3.3cm+\cos^4\theta
R_{\beta\bar{\beta}\alpha\bar{\alpha}}+\cos^3\theta\sin\theta
R_{\beta\bar{\beta}\alpha\bar{\beta}}+\cos^3\theta\sin\theta
R_{\beta\bar{\beta}\beta\bar{\alpha}}\\[4mm]
&&\hskip 3.3cm+\cos^2\theta\sin^2\theta
R_{\beta\bar{\beta}\beta\bar{\beta}}\\[4mm]
&&\hskip
2.9cm=\cos^2\theta\sin^2\theta(R_{\alpha\bar{\alpha}\alpha\bar{\alpha}}
+R_{\beta\bar{\beta}\beta\bar{\beta}}).
\end{array}$$
So we have
$R_{\beta\bar{\beta}\beta\bar{\beta}}+R_{\alpha\bar{\alpha}\alpha\bar{\alpha}}=0$,
if we choose $\theta$ such that $\cos^2\theta\sin^2\theta\neq 0$.
And this contradicts with the fact that $(M,g_{i\bar{j}}(t))$ has
positive holomorphic sectional curvature. Hence we proved that
$R_{\alpha\bar{\alpha}\beta\bar{\beta}}>0$, for all $t\in
(0,\delta)$.

Therefore by the Frankel conjecture, we know that $M$ is
biholomorphic to the complex projective space $CP^n$.

This completes the proof of Theorem 1.2.$$\eqno \#$$

\end{document}